\documentclass[12pt,leqno,twoside]{amsart}
\usepackage{amssymb,amsfonts,amsmath,amsthm}
\usepackage{latexsym}
\usepackage{verbatim}
\usepackage{graphicx}

\usepackage[poly,graph,tips]{xy}
\usepackage[hypertex]{hyperref}

\newtheorem{theorem}{Theorem}[section]

\newtheorem{lemma}[theorem]{Lemma}
\newtheorem{proposition}[theorem]{Proposition}
\theoremstyle{remark}
\newtheorem{remark}[theorem]{\sc Remark}
\theoremstyle{remark}

\theoremstyle{definition}
\newtheorem{definition}[theorem]{Definition}
\newtheorem{generic}[theorem]{Genericity condition}
\theoremstyle{remark}

\theoremstyle{remark}

\theoremstyle{remark}

\numberwithin{equation}{section}

\newcommand{\bR}{{\mathbb R}}

\newcommand{\bZ}{{\mathbb Z}}

\newcommand{\Tij}{P_iP_j}
\newcommand{\Tijk}{P_iP_jP_k}
\newcommand{\Pij}{P_{ij}}
\newcommand{\Pijk}{P_{ijk}}
\newcommand{\PM}{P_{1234}}

\newcommand{\PPPP}{P_1,P_2,P_3,P_4}
\newcommand{\TT}{{\mathbb T}}

\newcommand{\Vij}{V_{ij}}
\newcommand{\Vijk}{V_{ijk}}
\DeclareMathOperator{\Terr}{Terr}
\DeclareMathOperator{\CH}{CH}

\begin{document}

\title[The nine Morse generic tetrahedra]
 {The nine Morse generic tetrahedra}

\author{\sc D. Siersma}  

\address{Mathematisch Instituut, Universiteit Utrecht, PO
Box 80010, \ 3508 TA Utrecht
 The Netherlands.}

\email{siersma@math.uu.nl}

\author{\sc M. van Manen}

\address{Mathematical Institute, Hokkaido University, Kita 10, Nishi 8, Kita-Ku, Sapporo, Hokkaido, 060-0810, Japan}

\email{manen@math.sci.hokudai.ac.jp}

%\thanks{}

\subjclass[2000]{Primary 51M20; Secondary 68U05}

\keywords{Morse theory, polyhedra, Voronoi diagram, tetrahedral shape}

\date{\today}

%\dedicatory{}

%\commby{}

%%% ----------------------------------------------------------------------

\begin{abstract}
From computational geometry comes the notion of a Gabriel graph
of a point set in the plane. The Gabriel graph consists of those
edges connecting two points of the point set such that the circle
whose diameter is the edge does not contain any point of the point
set in its interior. We define a generalization
of the Gabriel graph to $n$ dimensions: the Morse poset.
Using Morse theory we prove that for a generic set of $4$
points in $\bR^3$ there are nine 
different Morse posets, up to combinatorial equivalence.
At the end we mention some open questions and report on the
results of computer experiments concerning these.  
%OPLOSSING: hier is dus een zin toegevoegd. 
We also compare our shape classification to another criterion
widely used in computer science.
\end{abstract}

%%% ---------------------------------------------------
\maketitle
%%% ---------------------------------------------------

%\tableofcontents
\setcounter{section}{0}

%%%%%%%%%%%%%%%%%%%%%%%%

\section{Introduction} 
Take $N$ points $P_1, \cdots , P_N$ in $\bR^n$ and consider 
the function 
$ d : \bR^n \to \bR $
defined by:
\begin{equation*}
d(X) = \min_{j=1 , \cdots , N} \ d(X,P_j)     
\end{equation*}
We want to study the evolution of the sets
$d_\epsilon = \{X \mid  d(X) \leq \epsilon\}$, as $\epsilon$ increases.
In particular we are interested in the Euler characteristic $\chi$ of $d_\epsilon$.
In case $\epsilon$ is very small $d_\epsilon$ consists of $N$ small
solid spheres. Thus $\chi=N$. If $\epsilon$ is very big, then $d_\epsilon$ is
contractible and hence $\chi=1$.
\par
For a generic set of points, $d$ is a topological Morse function. In that case,
as $\epsilon$ grows, $d$ passes through a number
of non-degenerate critical values.
When $d$ passes a critical value of index $i$, an $i$-cell gets attached.
\par
The number of critical points of index $i$ is $a_i$. From Morse theory we know
that 
\begin{equation*}
  \sum (-1)^ia_i = 1
\end{equation*}
As an example, take a triangle with an obtuse angle in the plane.
This is a special case of the above problem with $n=2$ and $N=3$. Assume further that
the two legs that encompass the obtuse angle have different lengths.
In that case $a_0=3$, $a_1=2$ and $a_2=0$.
For an acute triangle where the edges have different lengths,
we obtain $a_0=3$, $a_1=3$, and $a_2=1$.
In this sense, there are two different generic triangles.
\par
Returning to the $n$-dimensional case, each critical point of index
$i$ corresponds to a
subset of length $i$ of $\{P_1 , \cdots , P_N\}$, but not
every subset of length
$i$ corresponds to a critical point of index $i$.
( With the obtuse triangle the $2$-face of the triangle, that is:
the triangle itself, does not correspond to a local maximum
of the function $d$. )
The subsets of length $i$ that do correspond to critical points of $d$
will be called {\em active}.
An active subset defines a geometric $(i-1)$-simplex which we also
will call {\em active}.
Thus to the $N$-point set  $P_1 , \cdots , P_N$ we can associate a
set of subsets that 
are the active faces.
This set of subsets is partially ordered by inclusion, and thus
it is a poset. We will call it the {\em Morse poset}.
\par
{\em Question:} in $\bR^n$, for generic sets of $N$ points, how many different
Morse posets are there, up to combinatorial equivalence?
\par
When $N=n+1$ we thus ask how many ``different'' generic simplices there are.
This is a natural first problem to consider.
\par
In the plane the answer is two, see \cite{VoronoiDist} section 2.
For obtuse triangles, the Morse poset is
\begin{equation*}
  \{ \{ P_1 \},\{ P_2 \},\{ P_3 \},\{ P_1,P_3 \},\{ P_2,P_3 \}\}
\end{equation*}
and for acute triangles, we get the Morse poset 
\begin{equation*}
  \{ \{ P_1 \},\{ P_2 \},\{ P_3 \},\{ P_1,P_2 \}, \{ P_1,P_3 \},\{ P_2,P_3 \},\{ P_1,P_2,P_3 \}\}
\textrm{.}
\end{equation*}
\par
The main theorem of this article says that in $\bR^3$, there are nine different
generic tetrahedra ($3$-simplices). 
\par
In the first section, we recall the relevant Morse theory.
Then we establish some notation and state what 
Voronoi diagrams and Gabriel graphs are.
Next, we state and prove the main theorem.
% OPLOSSING: de oorspronkelijke zin was
% In the last section we report the results of numerical experiments
% about higher dimensional results and transitions.
% Maar we doen geen numerieke experimenten met hoger dimensionale
% zaken.
In the last section we discuss transitions in relation to the
configuration space of four points in $\bR^3$.
We also report on numerical experiments concerning volume
data of the different compartments
of the configuration space where the Morse poset is of a certain type. 
Finally we compare our classification to the classification by shape
types in \cite{MR2002k:65206}.
\section{Genericity conditions}
We focus here on $4$ points in $\bR^3$, but most of the notations
and definitions have straightforward extensions to the $(n,N)$ general case.
\par
We write $\Pij $  for the middle of the interval $\Tij $ ,
$\Pijk $  for the center of the circumscribed circle of
the triangle $\Tijk $,  and $\PM $  for the center of the
circumscribed sphere of the tetrahedron.
For the tetrahedron
itself, that is to say the convex hull of $\{ P_1 , P_2, P_3 , P_4 \}$,
we will use the notation $\TT$.
\par
We impose the following 
\begin{generic}\label{sec:gener-cond}
We require the set of points $\PPPP$ to be in general position, so that
the convex hull of $\PPPP$ is $3$-dimensional.
Moreover the points $\Pijk $  do not to lie on one of the 
edges of the triangle $\Tijk $  and also $\PM $  does not lie
in one of the planes of the triangles $\Tijk $.
\end{generic}
The condition \ref{sec:gener-cond} means that the function $d$ is not
too badly behaved. To express more carefully what that means we recall
the definition of a topological Morse function, see \cite{MR22:4071}.
Let $P\in\bR^n$. And let $f$ be a continuous real-valued function on $\bR^n$. 
\begin{definition}\label{sec:def-morse-fun}
$f$ is topologically regular at $P$ if there is some neighborhood $U$ of $P$ and
a homeomorphism $\phi\colon U\circlearrowleft$ such that one of the components
of $\phi$ is $f$.
The function $f$ has a critical point at $P\in\bR^n$ if $f$ is not topologically regular
at $P$. In that case, $P$ is called a non-degenerate critical point of 
index $i$ if 
there is a neighborhood $U$ of $P$ and a homeomorphism
$\phi\colon U\circlearrowleft$ such that 
\begin{equation*}
 f\circ\phi = f(P) - \sum_{j=1}^i x_j^2  + \sum_{j=i+1}^n x_j^2
\end{equation*}
A topological Morse function is a continuous function that has only
non-degenerate critical points.
\end{definition}
For topological Morse functions the two crucial statements that hold
in the differentiable case - the regular interval theorem and the attachment
of cells, see chapter 5 in \cite{Milnor1} - are true as well, as Morse
proves in \cite{MR47:9631}.
So if we show that $d$ is topologically regular, we
can apply those theorems, just as was done in \cite{VoronoiDist}.
\par
We will need the following notations:
\par
Let
\begin{equation*}
\Terr(P_i) = \{ X \in \bR ^3 | d(X,P_i) \le d(X,P_k) \textrm{for all k} \}
\textrm{,}
\end{equation*}
and $V_i = \Terr(P_i)$, $V_{ij}= V_i \cap V_j$ ($i$ different from $j$),
$V_{ijk} = V_i \cap V_j \cap V_k $ (all three i,j,k different),
$V_{1234} = V_1 \cap V_2 \cap V_3 \cap V_4$.
\begin{proposition}\label{sec:paktiv}
The function $d$ is a topological Morse function if the
condition \ref{sec:gener-cond} is fulfilled. In that case,
$d$ is topologically regular in all points of
of $\bR^3$, except in $\PPPP$, where $d$ has a minimum and
(perhaps) in the points $\Pij $, $\Pijk $ and $\PM $ .
Moreover $d$ has 
\begin{itemize}
\item[-]a minimum exactly in the points $\PPPP$
\item[-] a 1-saddle (saddle point of index $1$) in $\Pij$ iff.\ $\Pij = \Vij \cap \Tij$
\item[-] a 2-saddle (saddle point of index $2$) in $\Pijk$ iff.\ $\Pijk = \Vijk \cap \Tijk$
\item[-] a maximum in $\PM $ iff.\ $\PM = V_{1234} \cap \TT $, 
equivalently $\PM \in \TT$.
\end{itemize}
\end{proposition}
\begin{proof}
If $x=P_i$ then because the points lie in general position all $P_j$ with
$i\neq j$ lie at some positive distance from $x$, so $d$ has a minimum
there. If $x\in\Terr(P_i)$ but $x\neq P_i$ then $x$ is obviously topologically
regular.
\par
The function $d$ restricted to the interior $\Vij$ has a minimum if
$\Pij$ lies in that interior. This can only happen when $\Pij\neq\Pijk$, which
is assured by \ref{sec:gener-cond}.
\par
In the directions, orthogonal to $\Vij$, $d$ decreases so we see that $d$ has
a critical point of index $1$ at $P_{ij}$. 
\par
The function $d$ restricted to the interior of $\Vijk$ has a minimum at 
$\Pijk$ if $\Pijk$ lies in that interior. This can only be the case
when $\Pijk\neq\PM$. In the directions orthogonal to $\Vijk$ $d$ decreases,
so $d$ has a critical point of index $2$ at $\Pijk$.
\par
Finally if $\PM=V_{1234}$, $d$ obviously has a maximum there.
\end{proof}
\begin{remark}
Closer inspection might yield that \ref{sec:gener-cond} is actually not
necessary. See \cite{VoronoiDist} for the $2$-dimensional case.
\end{remark}
Throughout this article we will always assume that point sets satisfy
the genericity condition \ref{sec:gener-cond}.
\par
If $d$ has a critical point as described by the conditions
of the proposition \ref{sec:paktiv} 
then
we say that the corresponding center is {\em active}.
As explained in the introduction, an active center determines a subset of the points $P_1,P_2,P_3,P_4$,
which we call {\em active subset}. The {\em{}Morse poset} is the set of active
subsets.
Two sets of points in $\bR^3$ are called
{\em combinatorially equivalent} if there exists
a bijection, that sends the active subsets onto each other.
We want to give a classification with respect to 
this equivalence relation.
\section{Morse theoretic possibilities}
From proposition \ref{sec:paktiv} we know the maximal number the critical points 
of each type. Moreover the Euler characteristic should be $+1$.
This gives a priori the following $9$ possibilities:
\[
\begin{array}{|c|c|c|c|}
\hline
m & s1 & s2 & M \\
\hline
4 & 6  & 4  & 1 \\
4 & 5  & 3  & 1 \\
4 & 4  & 2  & 1 \\
4 & 3  & 1  & 1 \\
4 & 2  & 0  & 1 \\
4 & 6  & 3  & 0 \\
4 & 5  & 2  & 0 \\
4 & 4  & 1  & 0 \\
4 & 3  & 0  & 0 \\
\hline
\end{array}
\]
But not all possibilities will occur.
Since we start with $4$ points and the result should be a connected space,
we need at least $3$ saddle points of index $1$.
This rules out the possibility $(4,2,0,1)$.
\par
We remind the reader that the definitions that follow assume the
genericity condition \ref{sec:gener-cond}. For the more general definitions
see \cite{Mark} or \cite{MR2002k:65206}.
\par
The {\em Voronoi tesselation} of a point set $\{P_1 , \cdots , P_N\}$ in 
$\bR^n$ consists of the union of the sets $\Terr(P_i)$, together with
their natural combinatorial structure. Its $(n-1)$-skeleton consists
of the intersections $\Terr(P_i)\cap\Terr(P_j)$ and is usually
called 
the {\em Voronoi diagram}
The {\em Delaunay ``triangulation''} of that point set is the simplicial
complex dual to the Voronoi tesselation. Its $1$-skeleton consists
of the line segments 
\begin{equation*}
  \{ P_iP_j \quad\mid\quad \dim(\Terr(P_i)\cap\Terr(P_j)) = n-1 \quad\}
\end{equation*}
\par
The {\em Gabriel graph} ( see \cite{Mark}, where it is defined in 
$\bR^2$ ) of a point set in $\bR^n$ is formed by using the points as
vertices and placing an edge between two vertices exactly when the sphere,
whose diameter is given by that edge, does not include any point from the set.
The Gabriel graph is a subset of the Delaunay triangulation. 
Clearly the Gabriel graph is always connected.
\par
Our Morse poset is related to these much-used notions from
computational geometry as follows.
If in our case the edge between $P_i$ and $P_j$ is part of the
Gabriel graph then, the point $P_{ij}$ is a saddle point of $d$. 
The Morse poset is a generalization of the Gabriel 
graph. The Gabriel graph consists of the active subsets of length $2$. 
Two Morse posets can only be combinatorially equivalent if the
underlying Gabriel graphs are the same. 
\par
\section{Activity conditions}
\begin{figure}[htbp]
  \centering
  \includegraphics[height=9cm]{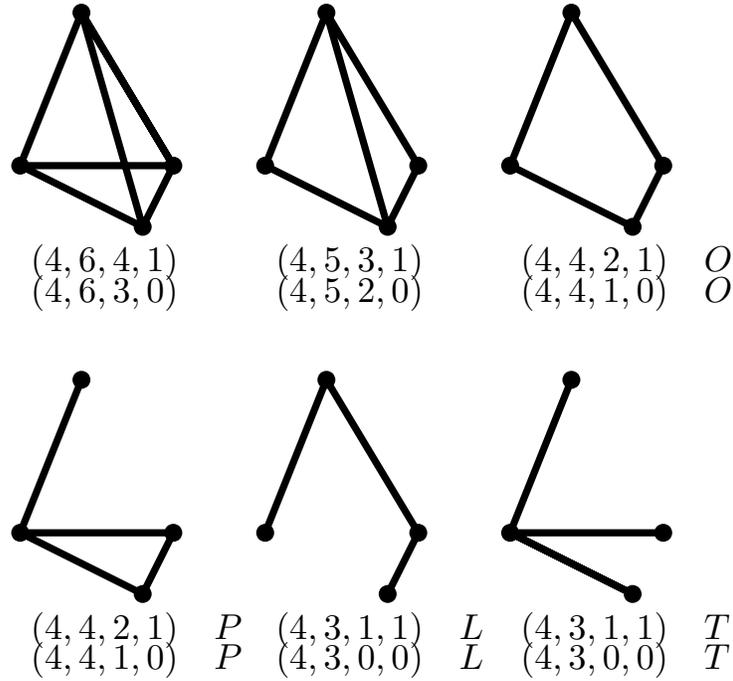}
  \caption{List of graphs}
  \label{fig:1}
\end{figure}

We list in figure \ref{fig:1} the (a priori) possible Gabriel
graphs for the above cases, they
are the connected graphs with $6$ vertices.
Just as the case $(4,2,0,1)$ can not occur, we will prove 
in this
section by deriving a contradiction that the case $(4,4,2,1)\, P$,
$(4,3,1,1)\, T$ and $(4,3,1,1)\, L$ do not occur. We will
prove our main theorem: 
\begin{theorem}\label{sec:main-thm}
Up to combinatorial equivalence of their Morse posets there are 
nine generic tetrahedra. They are uniquely described by the
nine Gabriel graphs $(4,6,4,1)$, $(4,6,3,0)$, $(4,5,3,1)$,
$(4,5,2,0)$, $(4,4,2,1)\, O$, $(4,4,1,0)\, O$, $(4,4,1,0)\, P$,
$(4,3,0,0)\, L$ and $(4,3,0,0)\, T$, drawn in figure
\ref{fig:1}. 
\end{theorem}
\begin{proof}
In the three cases that are to be excluded
$P_1P_2P_3P_4$ and $P_1P_2P_3$
are active. Let us see what these two conditions mean.
\subsection{Saddle points of index 1}
We consider the plane $E$ through $P_1P_2P_3$ 
(and put is as `ground plane' in the picture). 
The half space containing $P_4$ is called `above' , the other is called
`below'.\\
We are going to consider
the condition that the point $Y=P_{i4}$ is {\em active}.
We first make no assumptions about the triangle $P_1P_2P_3$,
except that it does not have a right angle.
We consider the point $X =P_4$ as a variable and denote its projection on the
$E$-plane by $X'$. Let $Y = P_{14}$ and its projection $Y'$.
From the condition that $Y = P_{14}$ is active we get:
\[
d(P_1,Y) = d(P_4,Y) \le \min{\{d(P_2,Y),d(P_3,Y)\}}
\]
For the projection, this means that in terms of Voronoi diagrams in $E$,
\[Y' \in \Terr(P_1) \]
In terms of the projection $X'$, this means:
\[X' \in 2 \Terr(P_1), \]
where with $2 \Terr(P_1)$ we mean scalar multiplication of $\Terr(P_1)$ 
by a factor $2$ from the point $P_1$.
\par
We can do the same for the activity of $P_{24}$ and $P_{34}$.
We get three subsets $2 \Terr(P_1)$, $2 \Terr(P_2)$, $2 \Terr(P_3)$, 
which cover the plane $E$. They divide the plane into regions,
where one, two or three of the points $P_{14}$, $P_{24}$ or $P_{34}$ 
are active.
\par
Inside the plane $E$, the border of $\Terr(P_i)$ consists of two
half-lines that meet in $P_{123}$. The scalar multiplication by
two maps $P_{123}$ to $P_i^*$, the antipodal point of $P_i$ 
on the circle through $P_1$, $P_2$ and $P_3$.
\par
Next, look at the case of an acute triangle drawn in figure \ref{fig:2}.
\par
The regions where only one $P_{i4}$ is active are both outside the
disc $D$, which is bounded by the circumscribed circle of triangle
$P_1P_2P_3$. 
\par
The picture for the obtuse case is in figure \ref{fig:3}.
\begin{figure}[thpb]
  \centering
  \includegraphics[height=9cm]{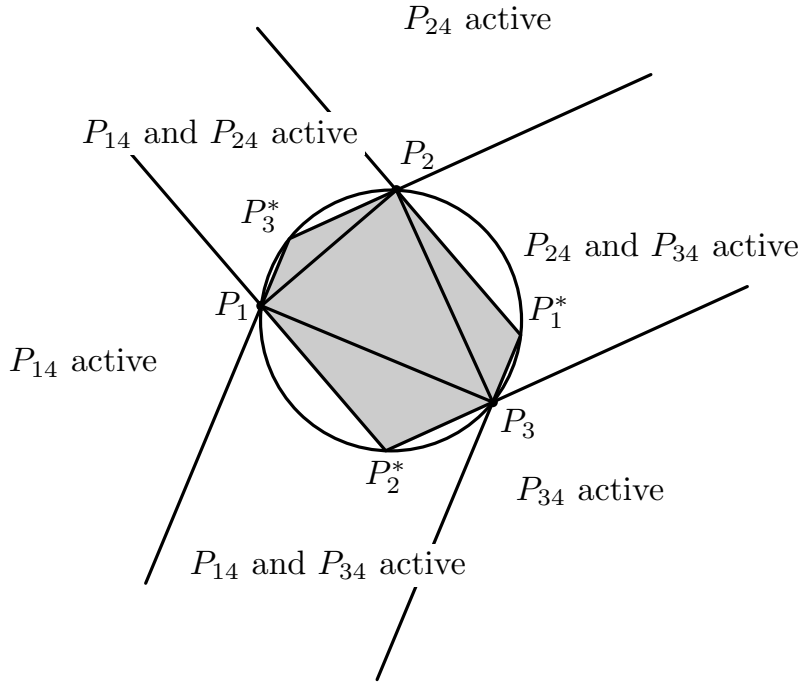}
  \caption{The plane $E$ when the triangle is acute}
  \label{fig:2}
\end{figure}
\begin{figure}[bhpt]
  \centering
  \includegraphics[height=8cm]{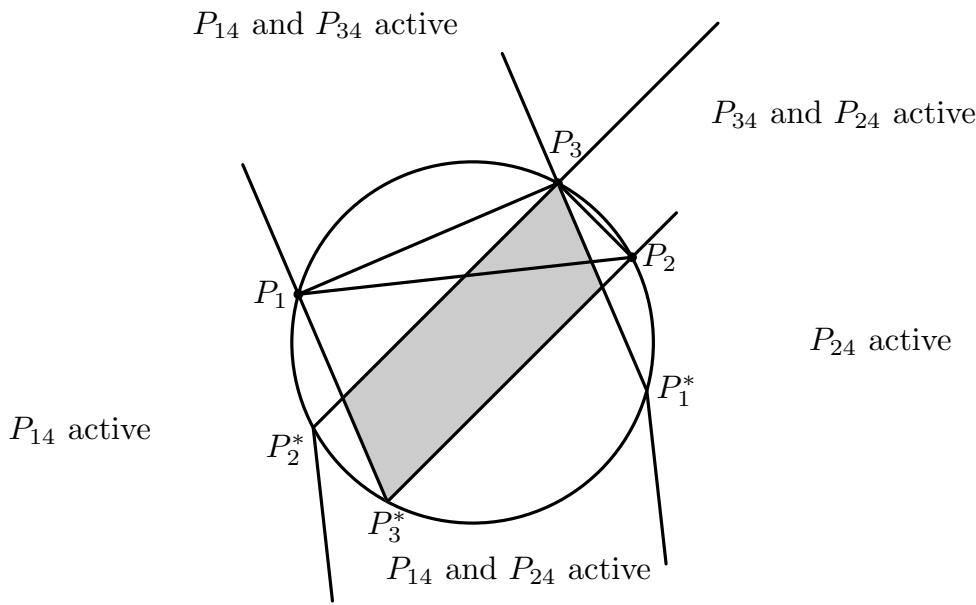}
  \caption{The plane $E$ when the triangle is obtuse}
  \label{fig:3}
\end{figure}
\par
Let $P_3$ the obtuse angle. We see that the region where only $P_{34}$
is active is outside $D$, but that the regions where $P_{14}$ or $P_{24}$
is the only active saddle point can have some intersections with $D$.
\subsection{The center $M$ of the circumscribed sphere.}
We assume again that $P_1P_2P_3$ lies in the plane $E$.
We know that $M=P_{1234}$ belongs to the axis of the triangle
$P_1P_2P_3$, so its
projection is $M^\prime = P_{123}$. Fix $M$ for the moment, and consider
$X=P_4$ as a variable.
Possible positions of $X=P_4$ are on the sphere with center
$M$ and radius
$r = d(P_1,M)=d(P_2,M)=d(P_3,M)$.
Denote by $P_1^*P_2^*P_3^*$, the mirror image of $P_1P_2P_3$ after central
reflection through $P_{123}$.
\par
In the cases we wish to exclude, $M$ is active, i.e. 
$M$ is a point, where the distance function has its maximum, i.e. $M$ lies
inside the tetrahedron. The point $X$ must lie in the cone with top $M$ and
base the triangle $P_1P_2P_3$, but also on the sphere $(M,r)$.
Consider the the plane $E"$,
which is the image of $E$ by central reflection through $M$. 
$E^{\prime\prime}$ is parallel to $E$. The cone intersects the plane in a triangle, 
which projects (orthogonally)
exactly onto $P_1^*P_2^*P_3^*$. $X$ belongs to a part of the cap of the
sphere, 
which lies ``above'' $E$ (seen from $M$). It follow that $X$ projects to a
point $X'$ which lies inside the disc $D$. 
(NB. the projection area contains at least the 
triangle $P_1^*P_2^*P_3^*$.) This is true for any position of $M$.
\subsection{Certain cases don't exist}
For the acute case this is all we need.
In the obtuse case though we can be more precise.
Let $P_3$ be the obtuse angle.
It is necessary for the activity of $M$ that $X'$ belongs to the circle sector
$P_2^* M' P_1^*$ (the one which does not contain $P_3$).
We combine the activity conditions for $M$ and one single $P_{i4}$.
\par
Start with an acute triangle in the ground plane. It follows that the
combination $M$ active and a single $P_{i4}$ active can not occur together.
\par
This rules out the graphs: 
$(4,3,1,1)\, T$ , $(4,3,1,1)\, L$ and $(4,4,2,1)\, P$.
\par
\subsection{Positions of index 2 saddles}
We are left with nine possibilities of Gabriel graphs.
The Gabriel graphs do not give a complete picture of the combinatorics of the
Morse points. The Gabriel graph, together with the information about activity of
$M$, tells us only the number of saddles of index 2, but not the position.
We discuss this now in detail and show that there are in all cases
unique positions for the saddles.
\par
$(4,6,4,1)$ : unique places (no choices).
\par
$(4,5,3,1)$ : There are two triangles, where all three midpoints 
are active.
Both triangles must be acute.
Take one of them in the `ground plane' and assume it to be $P_1P_2P_3$. 
Since $\PM$ is active, $\PM$ lies on
the same side as $P_4$. It follows that $V_{123}$ must intersect 
the ground plane
in the point $P_{123}$, so that point is active. 
The same reasoning applies to the second triangle; this fixes the place of two
2-saddles. The two positions left for the third saddle are 
combinatorially equivalent.
\par
$(4,4,1,0)\, O$ : all places equivalent.
\par
$(4,4,2,1)$ : Suppose $P_{12}$ is not active. Now $P_3$ or $P_4$ must 
be contained in $B=B(P_{12},r_{12})$, where $r_{12} = d(P_1,P_{12})$.
If $P_{123}$ is active, then triangle
$P_1P_2P_3$ must be acute, this means that
$P_3$ is outside the ball $B$. 
If $P_{124}$ is active, then triangle $P_1P_2P_4$ must be acute;
this means that
$P_4$ is outside the ball $B$. 
It follows that the situation where two active 2-saddles are ``separated''
by an non-active edge cannot occur.
This fixes the places of the 2-saddles up to 
permutation.
\par
$(4,6,3,0)$ : Since all midpoints of ribs are active,
we have that all triangles
are acute. Since  $\PM$ is not active it lies outside $\TT$.
There is exactly $1$ triangle such that $\PM$ and $P_4$ lie on different
sides of the plane of that triangle. The corresponding $P_{ijk}$ on that
axis cannot be active. This fixes the places of the other 2-saddles.
\par
$(4,5,2,0)$ : Choose a triangle, say $P_1P_2P_3$, where all three
edges are active. If $P_{123}$ is not active, then it follows that
$\PM$ lies below the ground plane $E$. We look again at the projection
$P_4'$. First $P_4'$ must lie inside $D$. But since two of the $P_{i4}$
(say $P_{14}$ and $P_{24}$) must be active we know that $P_4'$ also
must lie in the activity region, described above. The intersection is a subset
of $D$, which is contained in the region $D^*$ bounded by the 
arc $P_1P_3^*P_2$ and the interval $P_2P_1$.
There are still other conditions to meet: 
\begin{itemize}
\item[-] $P_4$ must lie
outside the ball $B(P_{12},r_{12})$, since $P_{12}$ is active.
\item[-] $P_4$ must lie inside the ball $B(P_{123},r_{123})$, 
where $r_{123}= d(P_1,P_{123})$.
\end{itemize}
This is not simultaneously possible if $P_4'$ lies in $D^*$.
It follows that there is only one possibility: the 2-saddles
are the centers of two triangles with all edges active.
\par
$(4,4,1,0)\, P$ : Let $P_1P_2P_3$ be the triangle with all three midpoints of
the edges active. This triangle is acute, take its plane as ground plane.
If $\PM$ is above the ground plane, then $P_{123}$ is active. Assume now 
$\PM$ is below the ground plane. It follows that $P_4 '$ lies inside the
disc $D$. But the fact that only $P_{14}$ is active means that $P_4 '$
is outside. This is a contradiction, so only $P_{123}$ can be active.
\par
$(4,3,0,0)\, T$ : all places equivalent. 
\par
$(4,3,0,0)\, L$ : all places equivalent. 
\par
That all these cases do occur follows from the computer experiments
described below. The proof of theorem \ref{sec:main-thm} is
complete.
\end{proof}
\section{Notes and remarks}
\subsection{Higher dimensional results}
We have not been able to prove a classification theorem in $\bR^n$.
However, upon request the authors will send interested readers a
computer program that calculates the list in higher dimensions 
by just trying a lot of random point sets; see below.
\par
Except for the results in \cite{VoronoiDist} on $4$ points in
the plane we have no results on the number of Morse posets for $N$ points
when $N>n+1$. 
\subsection{Relation with computer science criteria}
In the above it was remarked that the Morse poset is a
generalization of the Gabriel graph to higher dimensions.
Our notion of a Morse poset, being new, has not yet an
efficient stable algorithm to compute it for a general
point set. 
\subsection{The configuration space}
We consider the configuration space of $4$ points in $\bR^3$.
This is $\bR^3$ without the generalized diagonal
\begin{equation*}
  (\bR^3)^{(4)}
= \{ ( P_1, P_2 , P_3, P_4 ) \quad\mid\quad 
i\neq j \Rightarrow P_i\neq P_j \}
\end{equation*}
\begin{theorem}
The Morse poset does not change under scaling, rotation, or translation.
The quotient of $(\bR^3)^{(4)}$ by these group actions is a smooth
$5$-dimensional space.  
\end{theorem}
The proof is not very difficult, and we refrain from stating it here.
The non-generic tetrahedra form a hypersurface in
the configuration space. The theorem above says that there are
nine different types of compartments in the complement.
It does not give us any information on the number of components
of the complement and their topology.
\par
We do have the following information on the adjacency of the
types of compartments, for general $n$ and $N=n+1$.
%OPLOSSING: hier is de formulering dus nog voorzichtiger gemaakt.
\begin{theorem}\label{sec:configuration-space}
Let $\{P(\lambda)_i \}_{i=1, \cdots , n+1}$ be a generic smooth path in the
configuration space, that intersects the hypersurface exactly
once in the interval $[0,1]$. If the Morse poset changes there is a
$j$ such that 
the Morse poset of $\{P(1)_i \}_{i=1, \cdots , n+1}$ has one more ( resp.\
one less ) active subset of length $j$ and one more ( resp.\ one less ) active
subset $j+1$ than  $\{P(0)_i \}_{i=1, \cdots , n+1}$.
\end{theorem}
Again the proof is not very difficult, so we do not state it here.
\par
What the theorem says for instance is that the two types $(4,3,0,0)$ with
different Gabriel graphs are not adjacent.
\subsection{Statistics}
We carried out some statistical experiments to see how the different
types of tetrahedra are distributed
among different $4$-tuples of points. In this way complement
of the discriminant volume data are obtained.
Four points in $\bR^3$ do not form a bounded space, and thus
there lives no uniform probability distribution on it.
As the Morse poset is invariant under translations, scaling and rotations,
without loss of generality
we can assume that all $4$ points lie on $S^2$.
By translation and scaling, this can always be achieved. The
ratio of the ( infinite ) volumes of the complement of the discriminant space
will not change.
\par
A reference for the following material is \cite{MR82g:53003a}, chapter 9.
\par
In general, if one has a Riemannian manifold $(M,g)$ the metric
induces a
volume form
on $M$.
The simplest Riemannian metric comes when $M$ is embedded by a
map $\gamma$ as an
orientable hypersurface in $\bR^n$. Then the volume is given by:
\begin{equation*}
  \int_M \sqrt{\det( g_{ij} )} d x = \int_M \lVert \gamma_*e_1 \wedge
\cdots
\wedge \gamma_*e_{n-1} \rVert^{\frac12} d x
\end{equation*}
Suppose we have a map $ (\bR/\bZ)^{n-1} \rightarrow \bR^n$
that is onto $M\subset \bR^n$.
The uniform distribution on $[0,1]^{n-1}$ leads to a uniform
distribution on
 $M$ iff.
 $ \lVert \gamma_*e_1 \wedge \cdots \wedge \gamma_*e_{n-1}
\rVert^{\frac12}$
is a constant function on $M$.
%
%In general, if one has a Riemannian manifold $(M,g)$ the metric induces a volume form
%on $M$ by 
%\begin{equation}\label{eq:25}
%  \int_M \sqrt{\det( g_{ij} )} d x
%\end{equation}
%The simplest Riemannian metric comes when $M$ is embedded by a
%map $\gamma$ as an
%orientable hypersurface in $\bR^n$. Then the integral (\ref{eq:25}) reads
%\begin{equation*}
%  \int_M \lVert \gamma_*e_1 \wedge \cdots \wedge \gamma_*e_{n-1} \rVert^{\frac12} d x
%\end{equation*}
%Suppose we have a map $ (\bR/\bZ)^n \rightarrow \bR^n$
%that is onto $M\subset \bR^n$.
%The uniform distribution on $[0,1]^{n-1}$ leads to a uniform distribution on
% $M$ iff.\
% \begin{equation*}
%  \lVert \gamma_*e_1 \wedge \cdots \wedge \gamma_*e_{n-1} \rVert^{\frac12}
% \end{equation*}
%is a constant function on $M$.  
Let 
\begin{equation*}
\gamma(a_1,a_2) = ( \sin(\arccos(2a_1-1) )\sin(2\pi a_2) ,
\sin(\arccos(2a_1-1)) \cos(2 \pi a_2) , 2a_1  -1 )
\end{equation*}
With this choice it turns out that $d V = 4\pi d a_1 \wedge d a_2$.
So the induced probability measure on $S^2$ is uniform. 
We take two random numbers in $[0,1]$ and map them to $S^2$ using
$\gamma$.
\par
For our experiment we used the Gnu Scientific Library, see \cite{gsl}.
This library has an implementation of the apparently very reliable
MT19937 random number generator. 
We took samples of $10^8$ tetrahedra. Here is one.
\begin{center}
\begin{tabular}[hbtp]{|c|l|r|}
\hline
1 & $(4,3,0,0)\, L$ & 17,807,919 \\
2 & $(4,3,0,0)\, T$ &  898,689 \\
3 & $(4,4,1,0)\, O$ & 26,224,574 \\
4 & $(4,4,1,0)\, P$ & 16,421,773 \\
5 & $(4,5,2,0)$ &  24,350,101  \\
7 & $(4,4,2,1)$ & 3,266,345 \\
6 & $(4,6,3,0)$ & 1,797,721 \\
8 & $(4,5,3,1)$ & 2,697,783 \\
9 & $(4,6,4,1)$ & 6,535,095 \\
\hline
\end{tabular}
\end{center}
Other samples gave approximately the same numbers, with a maximum difference
of $3000$.
\par
For three random points on the circle, the chances are $50$ percent that one
gets an obtuse triangle. So these results are very different, and 
we have no explanation for them.
\subsection{Edelsbrunner ratio}
Denote, for brevity, $d_{ij}=d(P_i,P_j)$. The radius of the circumsphere
of $\TT$ is $R$.
\begin{definition}
The Edelsbrunner ratio $\rho$ is the circumradius $R$
divided by the minimal edge length $\min d_{ij}$.
\end{definition}
The ratio is used by Edelsbrunner ( see \cite{MR2002k:65206}, section 6.2 )
to classify
tetrahedra into ``shape types''. This article has the same objective,
so it is worthwhile to compare his criterion to ours. We will show
that small values of $\rho$ can be attained by all the types listed in 
theorem \ref{sec:main-thm}.
Depending on one's taste one can conclude 
that our classification is finer than the one proposed by
Edelsbrunner, or that it describes other features.
In any case, we hope to have improved on what
Edelsbrunner calls a ``fuzzy undertaking''.
\par
Up to rotations and translations,
the simplex is determined by the lengths of its six edges.
% OPLOSSING: zoals je ziet is de formulering hier nog iets
% cleaner geworden. 
We assume that the simplices satisfy the genericity condition
\ref{sec:gener-cond}.
The generic simplices are an open subset $V$ of $(\bR_{>0})^6$.
\begin{lemma}
Consider the circumradius $R$ as a function of the 
lengths $d_{ij}$. It is defined when the
four points are not coplanar.
It is homogeneous of degree $1$
and, where defined, it has nonzero gradient. In addition, on $V$, we have 
\begin{equation}\label{eq:1}
  \frac{\partial R}{\partial d_{ij}} > 0\quad 1 \leq i < j \leq 4
\end{equation}
\end{lemma}
\begin{proof}
The first two statements are obvious. For the gradient note that
\begin{equation*}
\sum_{1\leq i < j \leq 4} d_{ij}\frac{\partial R}{\partial d_{ij}} = 
 R(d_{12}, \cdots , d_{34})
\end{equation*}
because $R$ is homogeneous of degree one.
The circumradius is always $>0$, so the gradient is never zero. 
To see \eqref{eq:1}, fix the position of $P_1$, $P_2$ and $P_3$.
Varying only $d_{34}$ the point $P_4$ moves along a circle that
has axis $P_1P_2$. Along this circle, $R$ has a ``maximum'' when
the four points are coplanar, and a ``minimum'' when $P_3=P_4$.
Both cases correspond to non-generic simplices.
In between, the partial derivative of $R$ wrt.\ $d_{34}$
is clearly $>0$. 
\end{proof}
As a consequence, the Edelsbrunner ratio is a homogeneous
function of degree $0$. Hence, it is no restriction to 
assume $R=1$. The above lemma says that $\{R=1\}$ is a smooth
manifold $W$ in $V$. The following lemma can also be found in
\cite{MR2002k:65206}.
\begin{lemma}\label{sec:edelsbrunner-ratio-2}
$\rho$ has a unique minimum on $W$ when all the $d_{ij}$ are
equal.
\end{lemma}
\begin{proof}
Note first that it follows from~\eqref{eq:1} that the $d_{ij}$ are
regular functions in the sense of \ref{sec:def-morse-fun}.
We thus have six regular functions on a five dimensional space of 
which $\rho$ is the maximum:
\begin{equation*}
  \rho = \max \left( 
\frac1{d_{12}}, \cdots , \frac1{d_{34}} \right)
=\frac1{\min( d_{12}, \cdots , d_{34} )} 
\end{equation*}
Suppose that, on an open set $V^\prime$ in $\bR^N$, we have
$N+1$ regular functions
$f_0, \cdots , f_N$. Let $f$ be the maximum:
 $ f=\max(f_0, \cdots , f_N)$.
We assert that the function $f$ can only have a minimum
at $x\in V^\prime$  when $f(x)=f_i(x)$ , for
$i=0, \cdots , N$.
\par
Indeed because the functions are regular we can assume that they are all
linear functions. Our assertion is trivial in that case. 
\par
It follows that $\rho$ can only attain its minimum when all $d_{ij}$ are
equal. 
\end{proof}
\begin{lemma}\label{sec:edelsbrunner-ratio-1}
The minimal edge length is always the length of an active edge.
\end{lemma}
\begin{proof}
Suppose the minimal edge length were achieved by a non-active edge,
with endpoints $X$ and $Y$.
The edge is not active, so inside the circumsphere corresponding to that
edge there is a point that is not an endpoint of the edge. However the distance
from that point to either one of the endpoints $X$ and $Y$ is smaller than twice the circumradius.
\end{proof}
\begin{theorem}
On each of the compartments of the configuration space corresponding
to theorem \ref{sec:main-thm} $\rho$ is bounded from below by the 
values in the table below. The infimum corresponds to the quadruple of
points in the third column. These quadruples lie on the discriminant
hypersurface from theorem \ref{sec:configuration-space}, except for the
case $4641$, which corresponds to the global minimum of $\rho$.
\par
\begin{tabular}[hbtp]{|c|c|l|l|}
\hline
Type & Rho & Infimum \\
\hline
$(4,3,0,0)\, L$ & $\frac12\sqrt3 $ & (0,0,0),(1,0,0),(1,1,0),(1,1,1) \\
$(4,3,0,0)\, T$ & $\frac12\sqrt3 $ & (1,0,0) , (0,1,0) , (0,0,1), (0,0,0) \\
$(4,4,1,0)\, P$ &$\sqrt{\frac7{12}}$ & ($\frac12\sqrt{3}$,$-\frac12$, 0),($-\frac12\sqrt3$,$-\frac12$), (0,1,0),(0,1,$\sqrt3$) \\
$(4,4,1,0)\, O$ & $\frac12\sqrt2$  &  (0,0,0),(1,0,0),(1,1,0),(0,1,0)\\
$(4,5,2,0)$ &$\frac12\sqrt2$  & (0,0,0),(1,0,0),(1,1,0),(0,1,0)\\
$(4,6,3,0)$ &$\frac12\sqrt2$  &( $\cos\alpha_j$,$\sin\alpha_j$,0) $j=1,\cdots ,3$, (0,0,1) \\
$(4,4,2,1)$ &$\frac12\sqrt2$   &(1,0,0) , (0,1,0), (0,-1,0) , ($\cos\alpha$,0,$\sin\alpha$)\\
$(4,5,3,1)$ &$\frac12\sqrt2$   &(1,0,0) , (0,1,0), (0,-1,0) , (0,0,1)\\
$(4,6,4,1)$ &$\frac14\sqrt6$   &(0,1,0),($\frac12\sqrt3$,$-\frac12$,$0$),($-\frac12\sqrt3$,$-\frac12$,$0$),(0,0,$\sqrt2$)\\
\hline
\end{tabular}
\par
For the case $(4,6,3,0)$ the triangle $P_1P_2P_3$ is acute. 
For the case $(4,4,2,1)$ the angle $\alpha$ should satisfy
$-\frac\pi{2}\leq\alpha\leq 0$.
\end{theorem}
\begin{proof}
The case $(4,6,4,1)$ has been handled in the above. 
\par
Because the minimum is unique, the other cases only have an infimum for
$\rho$. These infima should thus correspond to quadruples on the discriminant
hypersurface. 
\par
The cases where $M=0$ are characterized by
$P_{1234} \notin \CH(P_1, \cdots , P_4)$. We may assume that $P_{1234}$ is the
origin and that $R=1$. It readily follows that all edges have length $>\sqrt2$.
For the cases $(4,4,1,0)\, O$, $(4,5,2,0)$ and $(4,6,3,0)$ the configurations in the above table
show that this is actually the infimum.
\par
If there is one triangle with an obtuse angle, we obtain in the same way
that $\rho>\frac12\sqrt2$. Thus, the configurations for $(4,4,2,1)$ and
 $(4,5,3,1)$
in the above table are actually infima. 
\par
We now need to study the first three cases.
Let us start with $(4,4,1,0)\, P$. 
We will assume that the triangle $P_1P_2P_3$ is active as well as 
the edges $P_1P_2$, $P_1P_3$ , $P_1P_4$ and $P_2P_3$. 
Maximizing the minimal length of the four edges on the Gabriel
graph results in four
edges of equal length. 
Take an equilateral triangle in the plane:
\begin{equation*}
P_1=(0,1,0) ~
P_2=  (-\frac12\sqrt3,-\frac12, 0) ~ P_3= ( \frac12\sqrt3 , -\frac12, 0)
\end{equation*}
This triangle has sides of length $\sqrt3$. 
Drawing the figure 
\ref{fig:2} in this case gives figure \ref{fig:5}. Thus the projection
of $P_4$ to the $P_1P_2P_3$ plane should lie in the gray area
above $P_1$ in figure \ref{fig:5}.
\par
Hence all cases with the projection of $P_4$ in the colored
region and $d(P_1,P_4)=\sqrt3$ give $(4,4,1,0)\, P$.
The infimum occurs when the circumradius is
minimal, and this is when $P_4=(0,1,\sqrt3)$. In that case
\begin{equation*}
\rho=  \sqrt{\frac7{12}}
\end{equation*}
\begin{figure}[hbtp]
  \centering
  \includegraphics[height=4cm]{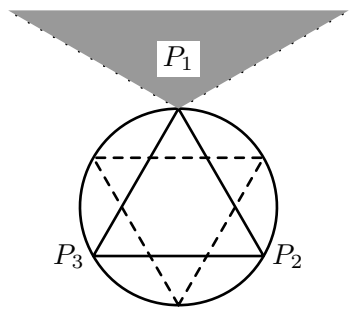}
  \caption{4410P Plane}
  \label{fig:5}
\end{figure}
The $(4,3,0,0)$ cases. 
There are two cases: either none of the triangles are acute, 
or at least one of the triangles is acute.
\par
If none of the triangles are acute we look at figure \ref{fig:3}. 
It follows that the projection of $P_4$ should lie outside of the triangle
 $P_1P_2P_3$. Hence the active edges are $P_1P_2$, $P_2P_3$ and 
$P_3P_4$. To maximize the minimal edge length of these
active edges, all three should have the same length, and the triangles
$P_1P_2P_3$, $P_2P_3P_4$ should be right angled. The 
simplex is then fixed by the values in the table above.
\par
In the second case, assume $P_1P_2P_3$ is acute. The triangle is
not active, so $P_4$ lies inside the sphere of $P_1P_2P_3$. Look at figure
\ref{fig:2}. 
The Gabriel graph contains the smallest edges, so all the other triangles
are obtuse angled. So the active edges are $P_1P_4$, $P_2P_4$ and $P_3P_4$.
To maximize the minimal length of these edges, they should all have equal
length. We see that the infimum for $\rho$ is as in the table above.
\end{proof}
\begin{remark} With a little heuristic reasoning, omitted here, it can be
argued that the configurations in the above table are actually unique infima.
\end{remark}
\bibliographystyle{amsalpha}
\bibliography{tetragen3}
\end{document}